\documentstyle[11pt,epic]{article}
\begin{document}

\newcommand{\ol}{\overline}
\newcommand{\mod}{{\rm mod}\ }
\newcommand{\W}{\hphantom{0}}
\newcommand{\bbZ}{{\makebox{\bf{Z}}}}
\newtheorem{lemma}{Lemma}[section]
\newtheorem{theorem}{Theorem}[section]
\newtheorem{corollary}{Corollary}[section]
\newcommand{\bbox}{\hspace*{6mm}\rule{2mm}{3mm}}
\newcommand{\I}{\Pi}
\newcommand{\QED}{{\hspace*{\fill}$\Box$}}

\renewcommand{\thefootnote}{\fnsymbol{footnote}}

\title{Some Integer Factorization Algorithms\\ using Elliptic Curves}

\author{Richard P. Brent\\
Computer Sciences Laboratory\\
Australian National University}

\date{24 September 1985\\%
Revised 10 Dec.~1985\\%
Republished 7 Nov.~1998%
\footnotetext{\hspace*{-18pt}Appeared in
{\em Australian Computer Science Communications\/} 8 (1986),
149--163.\\			
Retyped, with corrections and postscript, by
Frances Page at Oxford University Computing
Laboratory, 1998.
Copyright \copyright\ 1985, 1998 R.~P.~Brent.
\hspace*{\fill} rpb102 typeset using \LaTeX.}
}

\maketitle

\begin{abstract}
Lenstra's integer factorization algorithm is
asymptotically one of the fastest known algorithms, and is also ideally
suited for parallel computation.  We suggest a way in which the
algorithm can be speeded up by the addition of a second phase.  Under
some plausible assumptions, the speedup is of order $\log (p)$,
where $p$ is the factor which is found.  In practice the speedup
is significant.  We mention some refinements which give greater
speedup, an alternative way of implementing a second phase, and the
connection with Pollard's ``$p - 1$'' factorization algorithm.
\end{abstract}

\thispagestyle{empty}

\section{Introduction}

Recently H.W.~Lenstra Jr.~proposed a new integer factorization algorithm,
which we shall call ``Lenstra's algorithm'' or the ``one-phase elliptic
curve algorithm'' [17].  Under some plausible assumptions Lenstra's
algorithm finds a prime factor $p$ of a large composite integer $N$
in expected time

$$ T_1(p) = \exp\left(\sqrt{(2+o(1))\ln p \ln \ln p}\;\right),
\eqno{(1.1)}$$

\medskip \noindent
where ``$o(1)$'' means a term which tends to zero as $p \rightarrow \infty$.
Previously algorithms with running time
$\exp\left(\sqrt{(1 + o(1))\ln N \ln \ln N}\right)$
were known [27].  However, since $p^2 \le N$, Lenstra's
algorithm is comparable in the worst case and often much better, since it
often happens that $2\ln p \ll \ln N$.\\

The Brent-Pollard ``rho'' algorithm [5] is similar to Lenstra's algorithm
in that its expected running time depends on $p$, in fact it is of order
$p^{1/2}$.  Asymptotically $T_1(p) \ll p^{1/2}$, but because of the
overheads associated with Lenstra's algorithm we expect the ``rho''
algorithm to be faster if $p$ is sufficiently small.  The results of 
\S8 suggest how large $p$ has to be before Lenstra's algorithm
is faster.\\

After some preliminaries in \S\S2--4, we describe Lenstra's
algorithm in \S5, and outline the derivation of (1.1).  In
\S6 and \S7 we describe how Lenstra's algorithm can be speeded up
by the addition of a second phase which is based on the same idea as
the well-known ``paradox'' concerning the probability that two people
at a party have the same birthday [25].  The two-phase algorithm
has expected running time $O(T_1(p)/\ln p)$.  In practice, for $p$
around $10^{20}$, the ``birthday paradox algorithm'' is about 4 times faster 
than Lenstra's (one-phase) algorithm.  The performance of the various
algorithms is compared in \S8, and some refinements are mentioned
in \S9.

\section{Our unit of work}

The factorization algorithms which we consider use arithmetic
operations modulo $N$, where $N$ is the number to be factorized.  We
are interested in the case that $N$ is large (typically 50 to 200
decimal digits) so multiple-precision operations are involved.  As
our basic unit of work (or time) we take one multiplication modulo $N$
(often just called ``a multiplication'' below).  More precisely, given
integers $a$, $b$ in $[0, N)$, our unit of work is the cost of computing 
$a\!*\!b$ $\mod N$.  Because $N$ is assumed to be large, we can simplify the 
analysis by ignoring the cost of additions $\mod N$ or of
multiplications/divisions by small (i.e.~single-precision) integers,
so long as the total number of such operations is not much greater
than the number of multiplications $\mod N$.  See [13, 20] for
implementation hints.\\

In some of the algorithms considered below it is necessary to compute
inverses modulo $N$, i.e.~given an integer $a$ in $(0, N)$, compute
$u$ in $(0, N)$ such that $a\!*\!u = 1\; (\mod N)$.  We write 
$u = a^{-1}\; (\mod N)$.  $u$ can be computed by the extended GCD
algorithm [13], which finds integers $u$ and $v$ such that
$au + Nv = g$, where $g$ is the GCD of $a$ and $N$.  We can always
assume that $g = 1$, for otherwise $g$ is a nontrivial factor of $N$,
and the factorization algorithm can terminate.\\

Suppose that the computation of $a^{-1}\; (\mod N)$ by the extended
GCD algorithm takes the same time as $K$ multiplications $(\mod N)$.
Our first implementation gave $K \simeq 30$, but by using Lehmer's
method [16] this was reduced to $6 \le K \le 10$
(the precise value depending on the size of $N$).  It turns out that
most computations of $a^{-1}\; (\mod N)$ can be avoided at the
expense of about 8 multiplications $(\mod N)$, so we shall assume
that $K = 8$.\\

Some of the algorithms require the computation of large powers
$(\mod N)$, i.e.~given $a$ in $[0, n)$ and $b \gg 1$, we have to
compute $a^b \;(\mod N)$.  We shall assume that this is done by the
``binary'' algorithm [13] which requires between $\log_2b$ and
$2\log_2b$ multiplications $(\mod N)$~--~on average say
$(3/2)\log_2b$ multiplications (of which about $\log_2b$ are
squarings).  The constant 3/2 could be reduced slightly by use of
the ``power tree'' or other sophisticated powering algorithms~[13].

\section{Prime factors of random integers}

In order to predict the expected running time of Lenstra's algorithm
and our extensions of it, we need some results on the distribution of
prime factors of random integers.  Consider a random integer close to
$M$, with prime factors $n_1\ge n_2 \ge \ldots\;\;$.
For $\alpha \ge 1$, $\beta\ge 1$, define

$$\rho(\alpha) = \lim_{\scriptscriptstyle M \rightarrow \infty} {\rm Prob}\;
 \left(n_1 < M^{1/\alpha}\right)$$
\noindent
and
$$\mu(\alpha, \beta) = \lim_{\scriptscriptstyle M \rightarrow \infty} {\rm Prob}\;
\left(n_2 < M^{1/\alpha} \;\;{\rm and}\;\; n_1 < M^{\beta/\alpha}\right).$$

\medskip\noindent
(For a precise definition of ``a random integer close to $M$'', see [14].
It is sufficient to consider integers uniformly distributed in $[1, M].)$\\

Several authors have considered the function $\rho(\alpha)$, see for
example [7, 9, 13, 14, 18].  It satisfies a differential-difference
equation

$$\alpha\rho'(\alpha) + \rho(\alpha - 1) = 0$$

\medskip\noindent
and may be computed  by numerical integration from

$$\rho(\alpha) = \left\{
\begin{array}{ll}
1 &\; {\rm if}\;0\le\alpha \le1\\
\frac{1}{\alpha} \int^{\alpha}_{\alpha - 1} \rho(t)\;dt &\; {\rm if}\; 
\alpha > 1.
\end{array}
\right.$$

\medskip\noindent
We shall need the asymptotic results
$$ \ln \rho(\alpha) = -\alpha(\ln\alpha + \ln\ln\alpha - 1)
+ o(\alpha) \eqno{(3.1)} $$
\noindent
and
$$\rho(\alpha - 1)/\rho(\alpha) = \alpha(\ln\alpha
+ O(\ln\ln\alpha)) \eqno{(3.2)}$$

\medskip\noindent
as $\alpha \rightarrow \infty$.\\

The function $\mu(\alpha, \beta)$ is not so well-known, but is crucial
for the analysis of the two-phase algorithms.  Knuth and Trabb Pardo 
[14] consider $\mu(\alpha, 2)$ and by following their argument with 
trivial modifications we find that 

$$ \mu(\alpha, \beta) = \rho(\alpha)\, + \int^{\alpha - 1}_{\alpha - \beta}
\frac{\rho(t)}{\alpha - t}\;dt.
\eqno{(3.3)} $$

\medskip\noindent
When comparing the two-phase and one-phase algorithms the ratio
$\rho(\alpha)/\mu(\alpha, \beta)$ is of interest, and we shall need
the bound

$$ \rho(\alpha)/\mu(\alpha, \beta) =
   O\left(\ln\alpha(\alpha \ln\alpha)^{-\beta}\right)
\eqno{(3.4)}$$

\medskip\noindent
as $\alpha \rightarrow \infty$, for fixed $\beta > 1$.

\section{The group of an elliptic curve (mod \bf{\em p})}

In this section we consider operations $\mod p$ rather than $\mod N$,
and assume that $p$ is a prime and $p\ge5$.  When applying the results
of this section to factorization, $p$ is an (unknown) prime factor of
$N$, so we have to work $\mod N$ rather than $\mod p$.\\

Let $S$ be the set of points $(x, y)$ lying on the ``elliptic curve''

$$y^2 = x^3 + ax + b \hspace*{9mm} (\mod p),
\eqno{(4.1)}$$

\medskip\noindent
where $a$ and $b$ are constants, $4a^3 + 27b^2 \neq 0$. Let
$$ G = S\; \cup \; \{I\},$$

\medskip\noindent
where $I$ is the ``point at infinity'' and may be thought of as $(0, \infty)$.
Lenstra's algorithm is based on the fact that there is a natural way to
define an Abelian group on $G$.  Geometrically, if $P_1$ and 
$P_2 \in G$, we define $P_3 = P_1\!*\!P_2$ by taking $P_3$ to be
the reflection in the $x$-axis of the point $Q$ which lies on the
elliptic curve (4.1) and is collinear with $P_1$ and $P_2$.  Algebraically,
suppose $P_i = (x_i, y_i)$ for $i = 1, 2, 3$.  Then $P_3$ is defined by:

$$\begin{array}{l}
{\rm if}\; P_1 = I \;\;{\rm then}\; P_3 := P_2\\
{\rm else \ if}\; P_2 = I \;\;{\rm then}\; P_3 := P_1\\
{\rm else \ if}\; (x_1, y_1) = (x_2, -y_2) \;\;{\rm then}\; P_3 := I\\
{\rm else}\\
\hspace*{6mm}{\rm begin}\\
\hspace*{6mm}{\rm if}\; x_1 = x_2 \;\;{\rm then}\; \lambda 
              := (2y_1)^{-1}(3x^2_1 +a)\; \mod p\\
\hspace*{24mm}\;{\rm else}\; \lambda := (x_1 - x_2)^{-1}(y_1 - y_2)\; \mod p;\\
\hspace*{6mm}\{\lambda {\rm \ is\ the\ gradient\ of\ the\ line\ joining\ }
             P_1 {\rm\ and \ } P_2\}\\
\hspace*{6mm}x_3 := (\lambda^2 - x_1 - x_2)\; \mod p;\\
\hspace*{6mm}y_3 := (\lambda(x_1 - x_3) - y_1)\; \mod p\\
\hspace*{6mm}{\rm end}.
\end{array}$$

\medskip\noindent
It is well-known that $(G,*)$ forms an Abelian group with identity
element $I$.  Moreover, by the ``Riemann hypothesis for finite fields''
[12], the group order $g = |G|$ satisfies the inequality

$$|g - p - 1| \;<\; 2\sqrt p. \eqno{(4.2)}$$

\medskip
Lenstra's heuristic hypothesis is that, if $a$ and $b$ are chosen at
random, then $g$ will be essentially random in that the results of 
\S3 will apply with $M = p$.  Some results of Birch [3] suggest
its plausibility.  Nevertheless, the divisibility properties of $g$
are not quite what would be expected for a randomly chosen integer
near $p$, e.g.~the probability that $g$ is even is asymptotically
$2/3$ rather than $1/2$.  We shall accept Lenstra's hypothesis as we
have no other way to predict the performance of his algorithm.
Empirical results described in \S8 indicate that the algorithms
do perform roughly as predicted.\\

Note that the computation of $P_1\!*\!P_2$ requires $(3+K)$ units of
work if $P_1 \neq P_2$, and $(4+K)$ units of work if $P_1 = P_2$.
(Squaring is harder than multiplication!) If we represent $P_i$ as
$(x_i/z_i,\; y_i/z_i)$ then the algorithm given above for the computation
of $P_1\!*\!P_2$ can be modified to avoid GCD computations; assuming
that $z_1 = z_2$ (which can usually be ensured at the expense of 2
units of work), a squaring then requires 12 units and a nonsquaring 
multiplication requires 9 units of work.\\

The reader who is interested in learning more about the theory of
elliptic curves should consult [11], [12] or [15].

\section{Lenstra's algorithm}

The idea of Lenstra's algorithm is to perform a sequence of
pseudo-random trials, where each trial uses a randomly chosen elliptic
curve and has a nonzero probability of finding a factor of~$N$.  Let
$m$ and $m'$ be parameters whose choice will be discussed later.  
To perform a trial, first choose $P = (x, y)$ and $a$ at random.  
This defines an elliptic curve
$$y^2 = x^3 + ax + b \hspace*{9mm} (\mod N)
\eqno{(5.1)}$$

\medskip\noindent
(In practice it is sufficient for $a$ to be a single-precision random
integer, which reduces the cost of operations in $G$; also, there
is no need to check if GCD $(N, 4a^3 + 27b^2) \neq 1$ as this is
extremely unlikely unless $N$ has a small prime factor.)  Next compute
$Q = P^E$, where $E$ is a product of primes less than $m$,
$$ E = \prod_{p_i \;{\rm prime},\; p_i < m}\; {p_i}^{e_i}\;,$$
\noindent
where
$$e_i = \lfloor \ln(m')/\ln(p_i)\rfloor.$$

\medskip\noindent
Actually, $E$ is not computed.  Instead, $Q$ is computed by repeated
operations of the form $P := P^k$, where $k = {p_i}^{e_i}$ is a prime
power less than~$m'$, and the operations on $P$ are performed in the
group $G$ defined in \S4, with one important difference. The
difference is that, because a prime factor $p$ of $N$ is not known, all
arithmetic operations are performed modulo $N$ rather than modulo $p$.\\

Suppose initially that $m' = N$.  If we are lucky, all prime factors
of $g = |G|$ will be less than~$m$, so $g|E$ and $P^E = I$ in the
group $G$.  This will be detected because an attempt to compute 
$t^{-1}\; (\mod N)$ will fail because GCD $(N, t) > 1$.  In this case the
trial succeeds.  (It may, rarely, find the trivial factor $N$ if all prime
factors of $N$ are found simultaneously, but we neglect this
possibility.)\\

Making the heuristic assumption mentioned in \S4, and neglecting
the fact that the results of \S3 only apply in the limit as
$M ({\rm or\ } p) \rightarrow \infty$, the probability that a trial
succeeds in finding the prime factor $p$ of $N$ is just $\rho(\alpha)$,
where $\alpha = \ln(p)/\ln(m)$.\\

In practice we choose $m' = m$ rather than $m' = N$, because this
significantly reduces the cost of a trial without significantly 
reducing the probability of success.  Assuming $m' = m$, well-known
results on the distribution of primes [10] give $\ln(E) \sim m$, so the
work per trial is approximately $c_1m$, where
$c_1 = (\frac{11}{3} + K) \frac{3}{2 \ln 2}$.  Here $c_1$ is the product
of the average work required to perform a multiplication in $G$ times
the constant $\frac{3}{2 \ln 2}$ which arises from our use of the binary
algorithm for computing powers (see \S2).  Since $m = p^{1/\alpha}$,
the expected work to find $p$ is
$$ W_1(\alpha) \sim c_1p^{1/\alpha}/\rho(\alpha).
\eqno{(5.2)}$$

\medskip
To minimise $W_1(\alpha)$ we differentiate the right side of (5.2)
and set the result to zero, obtaining $\ln(p) = -\alpha^2\rho'(\alpha)/
\rho(\alpha)$, or (from the differential equation satisfied by $\rho$),
$$ \ln p = \frac{\alpha\rho(\alpha - 1)}{\rho(\alpha)}\;.
\eqno(5.3)$$

\medskip
In practice $p$ is not known in advance, so it is difficult to choose
$\alpha$ so that (5.3) is satisfied.  This point is discussed in
\S8.  For the moment assume that we know or guess an approximation
to $\log(p)$, and choose $\alpha$ so that (5.3) holds, at least
approximately.  From (3.2),
$$\ln p = \alpha^2(\ln\alpha + O(\ln \ln\alpha))\;,
\eqno(5.4)$$
\noindent
so
$$ \alpha \sim \sqrt{\frac{2 \ln p}{\ln \ln p}} \eqno{(5.5)}$$
\noindent
and
$$ \ln \; W_1(\alpha) \sim \frac{\rho(\alpha - 1)}{\rho(\alpha)}
-\ln\rho(\alpha) \sim 2 \alpha \ln \alpha
\sim \sqrt{2\ln p \ln\ln p}\;.
\eqno{(5.6)}$$
\noindent
Thus
$$T_1(p) = W_1(\alpha) =
\exp\left(\sqrt{(2+o(1))\ln p \ln \ln p}\right)\;,
\eqno{(5.7)}$$

\medskip\noindent
as stated in \S1.  It may be informative to write (5.7) as
$$T_1(p)  = W_1(\alpha) = p^{2/\alpha + o(1/\alpha)},
\eqno{(5.8)}$$

\medskip\noindent
so $2/\alpha$ is roughly the exponent which is to be compared with 1
for the method of trial division or $1/2$ for the Brent-Pollard 
``rho'' method.  For $10^{10} < p < 10^{30}$, $\alpha$ is in the
interval (3.2, 5.0).

\section{The ``birthday paradox'' two-phase algorithm}

In this section we show how to increase the probability of success
of a trial of Lenstra's algorithm by the addition of a ``second phase''.
Let $m = p^{1/\alpha}$ be as in \S5, and $m' = m^\beta >m$.
Let $g$ be the order of the random group $G$ for a trial of Lenstra's
algorithm, and suppose that $g$ has prime factors $n_1 \ge
n_2 \ge \ldots \;\;$ Then, making the same assumptions as in \S5,
the probability that $n_1 < m'$ and $ n_2 < m $ is $\mu(\alpha, \beta)$,
where $\mu$ is defined by (3.3).   Suppose we perform a trial of
Lenstra's algorithm, computing $Q = P^E$ as described in \S5.  With
probability $\mu(\alpha, \beta) - \rho(\alpha)$ we have $m \le
n_1 < m'$ and $n_2 < m$, in which case the trial fails because 
$Q \neq I$, but $Q^{n_1} = I$ in $G$.  (As in \S5, $g$ should
really be the order of $P$ in $G$ rather than the order of $G$, but
this difference is unimportant and will be neglected.)\\

Let $H = \langle Q \rangle$ be the cyclic group generated by $Q$.  A nice idea is
to take some pseudo-random function $f\!\!: Q \rightarrow Q$, define
$Q_0 = Q$ and $Q_{i+1} = f(Q_i)$ for $i = 0, 1, \ldots ,$ and generate
$Q_1, Q_2, \ldots$ until $Q_{2i} = Q_i$ in $G$.  As in the Brent-Pollard
``rho'' algorithm [5], we expect this to take $O(\sqrt n_1)$ steps.
The only flaw is that we do not know how to define a suitable
pseudo-random function~$f$.  Hence, we resort to the following (less
efficient) algorithm.\\

Define $Q_1 = Q$ and

$$Q_{j+1} = \left\{ \begin{array}{ll}
                     Q^2_j & {\rm with\ probability\ } 1/2,\\[.5ex]
                     Q^2_j\!*\!Q & {\rm with\ probability\ } 1/2,
                     \end{array}\right.
$$

\medskip\noindent
for $j = 1, 2, \ldots, r - 1$, so $Q_1, \ldots, Q_r$ are essentially
random points in $H$ and are generated at the expense of $O(r)$
group operations.  Suppose $Q_j = (x_j, y_j)$ and let

$$ d = \prod^{r-1}_{i=1}\; \prod^r_{j = i+1}\; (y_i - y_j)
\hspace*{9mm}(\mod N)
\eqno{(6.1)}$$

\medskip\noindent
If, for some $i<j\le r$, $Q_i=Q_j$ in $G$, then $p|(y_i - y_j)$
so $p|d$ and we can find the factor of $p$ of $N$ by computing GCD
$(N,d)$.  (We cannot find $i$ and $j$ by the algorithm used in the
Brent-Pollard ``rho'' algorithm because $Q_i = Q_j$ does not imply
that $Q_{i + 1} = Q_{j + 1}$.)\\

The probability that $p|d$ is the same as the probability that at
least two out of $r$ people have the same birthday (on a planet with
$n_1$ days in a year).  For example, if $n_1=365$ and $r=23$, the
probability $P_E \cong 1/2$.\\

\pagebreak[4]
In general, for $r \ll n_1$,

$$P_E = 1 - \prod^{r-1}_{j=1} (1 - j/n_1) \cong 1 - 
\exp\left(-\; \frac{r^2}{2n_1}\right),
\eqno{(6.2)}$$

\medskip\noindent
so we see that $P_E \ge 1/2$ if $r\; {\scriptstyle\stackrel{>}{\sim}}
\left(2n_1\ln2\right)^{1/2}$.\\

We can obtain a good approximation to the behaviour of the ``birthday
paradox'' algorithm by replacing the right side of (6.2) by a step
function which is 1 if $r^2 > 2n_1\ln2$ and 0 if $r^2 \le
2n_1\ln2$.  Thus, a trial of the ``birthday paradox'' algorithm will
succeed with probability approximately $\mu(\alpha, \beta)$, where
$\beta$ is defined by $r^2 = 2m^\beta\ln2$, i.e.

$$
\beta = \frac{2\ln r - \ln(2\ln2)}{\ln m} 
\eqno{(6.3)}$$

\medskip\noindent
and $\mu(\alpha, \beta)$ is as in \S3.  A more precise expression
for the probability of success is

$$\rho(\alpha) + \int^{\alpha-1}_0
\left\{1 - 2^{-p^{(t+\beta-\alpha)/\alpha}}\right\}
\frac{\rho(t)}{\alpha - t} \;\; dt.
\eqno{(6.4)}$$

\medskip\noindent
Computation shows that (6.3) gives an estimate of the probability of
success which is within 10\% of the estimate (6.4) for the
values of $p, \alpha$ and $\beta$ which are of interest, so we shall
use (6.3) below (but the numerical results given in \S8 were
computed using (6.4)).\\

A worthwhile refinement is to replace $d$ of (6.1) by

$${D} = \prod^{r-1}_{i=1}\; \prod^r_{j=i+1} \;(x_i - x_j)\hspace*{9mm} (\mod N).
\eqno{(6.5)}$$ 

\medskip\noindent
Since $(x_j, -y_j)$ is the inverse of $(x_j, y_j)$ in $H$, this
refinement effectively ``folds'' $H$ by identifying each point in $H$
with its inverse.  The effect is that (6.2) becomes

$$P_E \cong 1 - \exp\!\left(-\; \frac{r^2}{n_1}  \right)
\eqno{(6.6)}$$

\medskip\noindent
and (6.3) becomes $r^2 = m^\beta\ln2$, i.e.

$$
\beta = \frac{2\ln r - \ln\ln2}{\ln m}  
\eqno{(6.7)}$$

\medskip\noindent
(6.4) still holds so long as $\beta$ is defined by (6.7) instead of
(6.3).

\section{The use of fast polynomial evaluation}

Let $P(x)$ be the polynomial with roots $x_1, \ldots, x_r$, i.e.

$$ P(x) = \prod^r_{j=1} \;(x-x_j) = \sum^{r-1}_{j=0}\; a_jx^j 
\hspace*{9mm}(\mod N) \eqno{(7.1)}$$

{\samepage
\smallskip\noindent
and let $M(r)$ be the work necessary to multiply two polynomials of
degree $r$, obtaining a product of degree $2r$.  As usual, we assume
that all arithmetic operations are performed modulo $N$, where $N$ is
the number which we are trying to factorize.
}

Because a suitable root of unity $(\mod N)$ is not  known, we are
unable to use algorithms based on the FFT [1].  However, it is still
possible to reduce $M(r)$ below the obvious $O(r^2)$ bound.  For
example, binary splitting and the use of Karatsuba's idea [13] gives
$M(r) = O(r^{\log_2{3}})$.\\

The Toom-Cook algorithm [13] does not depend on the FFT, and it shows
that

$$M(r) = O\left(r^{1+(c/\ln r)^{1/2}}\right)
\eqno{(7.2)}$$

\medskip\noindent
as $r \rightarrow \infty$, for some positive constant $c$.  However,
the Toom-Cook algorithm is impractical, so let us just assume that we
use a polynomial multiplication algorithm which has

$$M(r) = O\left(r^{1+\varepsilon}\right)
\eqno{(7.3)}$$

\medskip\noindent
for some fixed $\varepsilon$ in $(0, 1)$.  Thus, using a recursive 
algorithm, we can evaluate the coefficients $a_0, \ldots, a_{r-1}$
of (7.1) in $O(M(r))$ multiplications, and it is then easy to obtain the
coefficients $b_j = (j+1)a_{j+1}$ in the formal derivative 
$P'(x) = \Sigma b_jx^j$.\\

Using fast polynomial evaluation techniques [4], we can now evaluate
$P'(x)$ at $r$ points in time $O(M(r))$.  However,

$$ {D}^2 = \prod^r_{j=1} \; P'(x_j),
\eqno{(7.4)} $$

\medskip\noindent
so we can evaluate ${D}^2$ and then GCD $(N, {D}^2)$.\\

Thus, we can perform the ``birthday paradox'' algorithm in time
$O(m) + O(r^{1 + \varepsilon})$ per trial, instead of $O(m) + O(r^2)$ if
(6.5) is evaluated in the obvious way.  To estimate the effect of this
improvement, choose $\alpha$ as in \S5  and $\beta = 2/(1 + \varepsilon)$
so that each phase of the ``birthday paradox'' algorithm takes about the
same time.  From (3.4) we have

$$\frac{\rho(\alpha)}{\mu(\alpha, \beta)} =
 O\left(\frac{\ln \alpha}{(\alpha \ln \alpha)^{2/(1+\varepsilon)}}\right)
= O\left(\frac{\ln\ln p}{(\ln p \ln\ln p)^{1/(1+\varepsilon)}}\right)\;.
\eqno{(7.5)}$$

\medskip\noindent
Thus, for any $\varepsilon' > \varepsilon$, we have a speedup of at least order
$(\ln p)^{1/(1+\varepsilon')}$ over Lenstra's algorithm.  If we use (7.2)
instead of (7.3) we obtain a speedup of order $\ln p$ in the same
way.\\

Unfortunately the constants involved in the ``$O$'' estimates make the
use of ``fast'' polynomial multiplication and evaluation techniques
of little value unless $r$ is quite large.  If $r$ is a power of 2
and binary splitting is used, so $\varepsilon = \log_{2}3 - 1 \cong 0.585$
above, we estimate that ${D}^2$ can be evaluated in $8r^{1+\varepsilon}
+ O(r)$ time units, compared to $r^2 /2 + O(r)$ for the obvious
algorithm.  Thus, the ``fast'' technique may actually be faster if
$r \ge 2^{10}$.  From the results of \S8, this occurs if
$p \ge 10^{22}$ (approximately).

\section{Optimal choice of parameters}

In Table 1 we give the results of a numerical calculation of the
expected work $W$ required to find a prime factor $p$ of a large
integer $N$, using four different algorithms:

\begin{enumerate}
\item The Brent-Pollard ``rho'' algorithm [5], which may be considered
as a benchmark.
\item Lenstra's one-phase elliptic curve algorithm, as described in
\S5.
\item Our ``birthday paradox'' two-phase algorithm, as described in
\S6, with $\varepsilon = 1$.
\item The ``birthday paradox'' algorithm with $\varepsilon = 0.585$, as described
in \S7, with $r$ restricted to be a power of 2.
\end{enumerate}

\begin{table}[th]
\begin{center}
~\\[-2ex]
\begin{tabular}{|rrrrr|}
\hline
&&&&\\[-2ex]
\multicolumn{1}{|c}{\hspace*{4mm} $\log_{10}p$\hspace*{2mm}}&
\multicolumn{1}{c}{\hspace*{4mm} Alg.~1\hspace*{2mm}}&
\multicolumn{1}{c}{\hspace*{4mm} Alg.~2\hspace*{2mm}}&
\multicolumn{1}{c}{\hspace*{4mm} Alg.~3\hspace*{2mm}}&
\multicolumn{1}{c|}{\hspace*{4mm} Alg.~4\hspace*{2mm}}\\[1ex]\hline

\W6\W\W\W&3.49\W &4.67\W &4.09\W &4.26\W\\
\W8\W\W\W&4.49\W &5.38\W &4.76\W &4.91\W\\
10\W\W\W &5.49\W &6.03\W &5.39\W &5.53\W\\
12\W\W\W &6.49\W &6.62\W &5.97\W &6.07\W\\
14\W\W\W &7.49\W &7.18\W &6.53\W &6.60\W\\
16\W\W\W &8.49\W &7.71\W &7.05\W &7.12\W\\
18\W\W\W &9.49\W &8.21\W &7.56\W &7.59\W\\
20\W\W\W &10.49\W&8.69\W &8.04\W &8.05\W\\
30\W\W\W &15.49\W&10.85\W&10.22\W&10.14\W\\
40\W\W\W &20.49\W&12.74\W&12.11\W&11.97\W\\
50\W\W\W &25.49\W&14.44\W&13.82\W&13.62\W\\[1ex]\hline
\multicolumn{5}{|c|}{~}\\[-2ex]
\multicolumn{5}{|c|}{
Table 1: $\log_{10}W$ versus $\log_{10}p$ for Algorithms 1--4}\\[1ex]\hline
\end{tabular}
\end{center}
\end{table}

In all cases $W$ is measured in terms of multiplications $(\mod N)$,
with one extended GCD computation counting as 8 multiplications
(see \S2).  The parameters $\alpha$ and $\beta$ were chosen to
minimize the expected value of $W$ for each algorithm (using numerical
minimization if necessary).  The results are illustrated in Figure 1.\\

From  Table 1 we see that Algorithm 3 is better than Algorithm 1
for $p \;{\scriptstyle\stackrel{>}{\sim}}\; 10^{10}$, while Algorithm 2
is better than Algorithm 1 for $p \;{\scriptstyle\stackrel{>}{\sim}}\; 10^{13}$.
Algorithm 3 is 4 to 4.5 times faster than Algorithm 2.  Algorithm 4 is
slightly faster than Algorithm 3 if
$p \;{\scriptstyle\stackrel{>}{\sim}}\; 10^{22}$.\\

The differences between the algorithms appear more marked if we consider
how large a factor $p$ we can expect to find in a given time.  Suppose
that we can devote $10^{10}$ units of work to the factorization.
Then, by interpolation in Table 1 (or from Figure 1), we see that
the upper bounds on $p$ for Algorithms 1, 2 and 3 are about $10^{19}$,
$10^{26}$ and $10^{29}$ respectively.

\begin{table}[h]
\begin{center}
\begin{tabular}{|rrrrrrrrr|}
\hline
&&&&&&&&\\[-2ex]
\multicolumn{1}{|c}{$\log_{10}p$}&
\multicolumn{1}{c}{$\alpha$}&
\multicolumn{1}{c}{$\beta$}&
\multicolumn{1}{c}{$m$}&
\multicolumn{1}{c}{$r$}&
\multicolumn{1}{c}{$T$}&
\multicolumn{1}{c}{$w_{21}$}&
\multicolumn{1}{c}{$m/T$}&
\multicolumn{1}{c|}{$S$}\\[1ex] \hline
&&&&&&&&\\[-2ex]
10\W&3.72&1.56&484   & 104&12.1 &0.64& 40&4.37\\
20\W&4.65&1.35&19970 & 669&147.5&0.47&135&4.46\\
30\W&5.36&1.27&397600&2939&1141&0.44&348&4.32\\[1ex]\hline
\multicolumn{9}{|c|}{~}\\[-2ex]
\multicolumn{9}{|c|}
{Table 2: Optimal parameters for Algorithm 3}\\[1ex]\hline

\end{tabular}
\end{center}
\end{table}

{\samepage
In Table 2 we give the optimal parameters $\alpha$, $\beta$,
$m = p^{1/\alpha}$,
$r = (m^\beta \ln2)^{1/2}$,\\
$T =$ expected number of trials (from (6.4)),
$m/T$, $w_{21} =$ (work for phase 2)/(work for phase 1),\\
and $S =$ speedup over Lenstra's algorithm, all for Algorithm 3 and
several values of $p$.
}

\begin{table}[ht]
\begin{center}
\begin{tabular}{|cccc|}
\hline
&&&\\[-2ex]
Algorithm & Number of & Observed work & Expected work \\
         & Factorizations&per factor/$10^6$&per factor/$10^6$ \\[1ex]\hline
&&& \\[-2ex]
2&126& 3.41 $\pm$ 0.30 & 4.17 \\
3&100& 0.74 $\pm$ 0.06 & 0.94 \\[1ex] \hline
\multicolumn{4}{|c|}{~}\\[-2ex]
\multicolumn{4}{|c|}{Table 3: Observed versus expected work per factor
for Algorithms 2--3}\\[1ex]\hline
\end{tabular}
\end{center}
\end{table}

In order to check that the algorithms perform as predicted, we factored
several large $N$ with smallest prime factor $p \cong 10^{12}$.
In Table 3 we give the observed and expected work to find each factor
by Algorithms 2 and 3.  The agreement is reasonably good, considering the
number of approximations made in the analysis.  If anything the
algorithms appear to perform slightly better than expected.\\[2ex]

In practice we do not know $p$ in advance, so it is difficult to
choose the optimal parameters $\alpha, \beta$ etc.  There are several
approaches to this problem.  If we are willing to devote a certain
amount of time to the attempt to factorize $N$, and intend to give up
if we are unsuccessful after the given amount of time, then we may
estimate how large a factor $p$ we are likely to find (using Table 1
or Figure 1) and then choose the optimal parameters for this ``worst
case'' $p$.  Another approach is to start with a small value of $m$
and increase $m$ as the number of trials $T$ increases.  From Table 2,
it is reasonable to take $m/T \cong 135$ if we expect to find a prime
factor $p \cong 10^{20}$.  Once $m$ has been chosen, we may choose
$r$ (for Algorithms 3 or 4) so that $w_{21}$ (the ratio of the work
for phase 2 to the work for phase 1) has a moderate value.  From Table 2,
$w_{21} \cong 0.5$ is reasonable.  In practice these ``ad hoc''
strategies work well because the total work required by the algorithms
is not very sensitive to the choice of their parameters (e.g.~if
$m$ is chosen too small then $T$ will be larger than expected, but
the product $mT$ is relatively insensitive to the choice of $m$).

\begin{figure}[ht]
\begin{center}					
\setlength{\unitlength}{1mm}
\begin{picture}(120,100)(0,0)
\put(0,0){\line(1,0){120}}
\multiput(0,-2)(20,0){7}{\line(0,1){2}}
\put(-2.5,-6){\makebox(5,3){6}}
\put(17.5,-6){\makebox(5,3){8}}
\put(37.5,-6){\makebox(5,3){10}}
\put(57.5,-6){\makebox(5,3){12}}
\put(77.5,-6){\makebox(5,3){14}}
\put(97.5,-6){\makebox(5,3){16}}
\put(117.5,-6){\makebox(5,3){18}}
\put(0,0){\line(0,1){100}}
\multiput(-2,0)(0,10){11}{\line(1,0){2}}
\put(-7,-2.5){\makebox(5,3){0}}
\put(-7,7.5){\makebox(5,3){1}}
\put(-7,17.5){\makebox(5,3){2}}
\put(-7,27.5){\makebox(5,3){3}}
\put(-7,37.5){\makebox(5,3){4}}
\put(-7,47.5){\makebox(5,3){5}}
\put(-7,57.5){\makebox(5,3){6}}
\put(-7,67.5){\makebox(5,3){7}}
\put(-7,77.5){\makebox(5,3){8}}
\put(-7,87.5){\makebox(5,3){9}}
\put(-7,97.5){\makebox(5,3){10}}
\drawline(0,34.9)(120,94.9)
\put(120,94.9){\makebox(5,3){1}}
\drawline(0,46.7)(20,53.8)(40,60.3)(60,66.2)(80,71.8)(100,77.1)(120,82.1)
\put(120,82.1){\makebox(5,3){2}}
\drawline(0,40.9)(20,47.6)(40,53.9)(60,59.7)(80,65.3)(100,70.5)(120,75.6)
\put(120,75.6){\makebox(5,3){3}}
\end{picture}

\vspace*{10mm}

\begin{tabular}{ll}
\underline{Figure 1:} & $\log_{10}W$ versus $\log_{10}p$ for Algorithms 1--3\\
\end{tabular}
\end{center}					
\end{figure}

\section{Further refinements}

In this section we mention some further refinements which can be used
to speed up the algorithms described in \S\S5--7.  Details will
appear elsewhere.

\subsection{Better choice of random points}

Let $e \geq 1$ be a fixed exponent, let $b_i$ and $\overline{b}_i$
be random linear functions of $i$, $a_i = b^e_i$,
$\ol{a}_i = \ol{b}\:\!^e_i$,
and $rs \sim m^\beta$.  In the birthday paradox algorithm we may
compute

$$(x_i, y_i) = Q ^{a_i} \hspace*{24mm} (i = 1, \ldots, r)$$

\medskip\noindent
and

$$(\ol{x}_j, \ol{y}_j) = Q ^{\ol{a}_j} \hspace*{24mm} (j = 1, \ldots, s)$$

\medskip\noindent
and replace (6.1) by

$$ d = \prod^s_{j=1} \; \prod^r_{i=1}\; (x_i - \ol{x}_j)
\hspace*{9mm}(\mod N)
\eqno{(9.1)}$$

\medskip\noindent
Using $e > 1$ is beneficial because the number of solutions of
$x^e = 1\; (\mod n_1)$ is GCD $(e, n_1  - 1)$.  We take $b_i$ and
$\ol{b}_i$ to be linear functions of $i$ so that the $e$-th
differences of the $a_i$ and $\ol{a}_i$ are constant, which allows 
the computation of $x_1, \ldots, x_r$ and  $\ol{x}_1, \ldots, \ol{x}_s$
in $O((r+s)e)$ group operations.  The values of $\ol{x}_j$ do not
need to be stored, so storage requirements are $O(r)$ even if
$s \gg r$.  Moreover, by use of rational preconditioning [22, 29]
it is easy to evaluate (9.1) in $(r+O(\log r))s/2$ multiplications.
Using these ideas we obtain a speedup of about 6.6 over the
one-phase algorithm for $p \cong 10^{20}$.

\subsection{Other second phases}

Our birthday paradox idea can be used as a second phase for Pollard's
``$p - 1$'' algorithm [23].  The only change is that we work over a
different group.  Conversely, the conventional second phases for
Pollard's ``$p - 1$'' algorithm can be adapted to give second phases
for elliptic curve algorithms, and various tricks can be used to speed
them up [19].  Theoretically these algorithms give a speedup of the
order $\log \log (p)$ over the one-phase algorithms, which is not
as good as the $\log (p)$ speedup for the birthday paradox algorithm [6].
However, in practice, the speedups are comparable (in the range 6 to 8).
We prefer the birthday paradox algorithm because it does not require
a large table
(or on-line generation)	
of primes for the second phase, so it is easier to
program and has lower storage requirements.

\subsection{Better choice of random elliptic curves}

Montgomery [21] and Suyama [28] have shown that it is possible to
choose ``random'' elliptic curves so that $g$ is divisible by certain
powers of 2 and/or 3.  For example, we have implemented a suggestion of
Suyama which ensures that $g$ is divisible by 12.  This effectively
reduces $p$ to $p/12$ in the analysis above, so gives a speedup which
is very significant in practice, although not significant 
asymptotically.

\subsection{Faster group operations}

Montgomery [21] and Chudnovsky and Chudnovsky [8] have shown that the
Weierstrass normal form (5.1) may not be optimal if we are interested in
minimizing the number of arithmetic operations required to perform
group operations.  If (5.1) is replaced by

$$by^2 = x^3 + ax^2 + x \hspace*{9mm}(\mod N)  
\eqno{(9.2)}$$

\medskip\noindent
then we can dispense with the $y$ coordinate and compute $P^n$ in
$10 \log_2n + O(1)$ multiplications $(\mod N)$, instead of about
$\frac{3}{2}(\frac{11}{3} + K)\log_2n$ multiplications $(\mod N)$, a saving of about
43\% if $K = 8$.

The effect of the improvements described in
\S\S9.3--9.4 is to speed up both the one-phase and two-phase
algorithms by a factor of 3 to 4.

\section{Conclusion}

Lenstra's algorithm is currently the fastest known factorization
algorithm for large $N$ having a factor
$p\; {\scriptstyle\stackrel{>}{\sim}}\; 10^{13}$,
$\ln p/\ln N \ll 1/2$.  It is also ideally suited to
parallel computation, since the factorization process involves a
number of independent trials which can be performed in parallel.\\

We have described how to improve on Lenstra's algorithm by the addition
of a second phase.  The theoretical speedup is of order $\ln(p)$.
From an asymptotic point of view this is not very impressive, but in
practice it is certainly worth having and may increase the size of
factors which can be found in a reasonable time by several orders of
magnitude (see Figure 1 and the comments in \S8).\\

Given increasing circuit speeds and increasing use of parallelism,
it is reasonable to predict that $10^{14}$ multiplications might be
devoted to factorizing a number in the not-too-far-distant future
(there are about $3 \times 10^{13}$ microseconds in a year).  Thus,
from Table 1, it will be feasible to find prime factors $p$ with up to
about 50 decimal digits by the algorithms based on elliptic curves.
Other algorithms [27] may be even more effective on numbers which are
the product of two roughly equal primes.  This implies that the composite
numbers $N$ on which the RSA public-key cryptosystem [25, 26] is based
should have at least 100 decimal digits if the cryptosystem is to be
reasonably secure.

\section{Acknowledgements}

I wish to thank Sam Wagstaff, Jr.~for introducing me to Lenstra's
algorithm, and Brendan McKay, Andrew Odlyzko, John Pollard, Mike
Robson and Hiromi Suyama for their helpful comments on a first draft
of this paper.

\section{References}

{%
\begin{enumerate}

\item A.~V.~Aho, J.~E.~Hopcroft and J.~D.~Ullman,
{\em The Design and Analysis of Computer Algorithms\/},
Addison-Wesley, 1974.

\item E.~Bach,
{\em Lenstra's Algorithm for Factoring with Elliptic Curves (expos\'{e})\/},
Computer\linebreak
Science Dept., Univ.~of Wisconsin, Madison, Feb.~1985.

\item B.~J.~Birch,
How the number of points of an elliptic curve over a fixed prime field
varies,
{\em J.~London Math.~Soc.\/} 43 (1968), 57--60.

\item A.~Borodin and I.~Munro,
{\em The Computational Complexity of Algebraic and Numeric\linebreak
Problems\/},
Elsevier, 1975.

\item R.~P.~Brent,
An improved Monte Carlo factorization algorithm,
{\em BIT\/} 20 (1980), 176--184.

\item R.~P.~Brent,
{\em Some integer factorization algorithms using elliptic curves\/},
Report CMA-R32-85, Centre for Math.~Analysis, Australian National
University, Sept.~1985, \S6.

\item N.~G.~de Bruijn,
The asymptotic behaviour of a function occurring in the theory of primes,
{\em J.~Indian Math.~Soc\/}. 15 (1951), 25--32.

\item D.~V.~Chudnovsky and G.~V.~Chudnovsky,
{\em Sequences of numbers generated by addition in formal groups and new
primality and factorization tests\/},
preprint, Dept.~of Mathematics, Columbia Univ., July 1985.

\item K.~Dickman,
On the frequency of numbers containing prime factors of a certain relative
magnitude,
{\em Ark.~Mat., Astronomi och Fysik\/},
22A, 10 (1930), 1--14.

\item G.~H.~Hardy and E.~M.~Wright,
{\em An Introduction to the Theory of Numbers\/},
Oxford\linebreak
University Press, 4th Edition, 1960.

\item K.~F.~Ireland and M.~Rosen,
{\em A Classical Introduction to Modern Number Theory\/},
Springer-Verlag, 1982, Ch.~18.

\item J-R.~Joly,
Equations et vari\'{e}t\'{e}s alg\'{e}briques sur un corps fini,
{\em L'Enseignement Math\'{e}matique\/} 19 (1973), 1--117.

\item D.~E.~Knuth,
{\em The Art of Computer Programming\/},
Vol.~2 (2nd Edition), Addison-Wesley, 1982.

\item D.~E.~Knuth and L.~Trabb Pardo,
Analysis of a simple factorization algorithm,
{\em Theoretical Computer Science\/} 3 (1976), 321--348.

\item S.~Lang,
{\em Elliptic Curves -- Diophantine Analysis\/},
Springer-Verlag, 1978.

\item D.~H.~Lehmer,
Euclid's algorithm for large numbers,
{\em Amer.~Math.~Monthly\/} 45 (1938), 227-233.

\item H.~W.~Lenstra, Jr.,
{\em Elliptic Curve Factorization\/},
personal communication {\em via\/} Samuel\linebreak
Wagstaff Jr., Feb.~1985.

\item J.~van de Lune and E.~Wattel,
On the numerical solution of a differential-difference equation
arising in analytic number theory,
{\em Math.~Comp.} 23 (1969), 417--421.

\item P.~L.~Montgomery,
{\em Speeding the Pollard methods of factorization\/},
preprint, System\linebreak
Development Corp., Santa Monica, Dec.~1983.

\item P.~L.~Montgomery,
Modular multiplication without trial division,
{\em Math.~Comp.} 44 (1985), 519--521.

\item P.~L.~Montgomery,
personal communication, September 1985.

\item M.~Paterson and L.~Stockmeyer,
On the number of nonscalar multiplications necessary to evaluate
polynomials,
{\em SIAM J.~Computing\/} 2 (1973), 60--66.

\item J.~M.~Pollard,
Theorems in factorization and primality testing,
{\em Proc.~Cambridge\linebreak
Philos.~Soc.\/} 76 (1974), 521--528.

\item J.~M.~Pollard,
A Monte Carlo method for factorization,
{\em BIT\/} 15 (1975), 331--334.

\item H.~Riesel,
{\em Prime numbers and computer methods for factorization\/},
Birkhauser, 1985.

\item R.~L.~Rivest, A.~Shamir and L.~Adleman,
A method for obtaining digital signatures and public-key cryptosystems,
{\em Comm.~ACM\/} 21 (1978), 120--126.

\item R.~D.~Silverman,
{\em The multiple polynomial quadratic sieve\/},
preprint, Mitre Corp., Bedford Mass., 1985.

\item H.~Suyama,
{\em Informal preliminary report {\rm(8)}\/,}
personal communication, October 1985.

\item S.~Winograd,
Evaluating polynomials using rational auxiliary functions,
{\em IBM Technical Disclosure Bulletin\/}
13 (1970), 1133--1135.

\end{enumerate}
}

\pagebreak[4]
\section*{Postscript and historical note (added 7 November 1998)}

The \LaTeX\ source file was retyped in 1998 from the version (rpb102)
which appeared in {\em Proceedings of the Ninth
Australian Computer Science Conference\/}, special issue of
{\em Australian Computer 
Science Communications\/} 8 (1986),
149--163 [submitted 24~September 1985, and in final form 10~December 1985].
No attempt has been made to update the contents, but minor
typographical errors have been corrected (for example, in equations
(1.1), (6.3), (6.7) and (9.2)). Some minor changes have been
made for cosmetic reasons, e.g. $d'$ was changed to $D$ in (6.5),
and some equations have been displayed more clearly using the
\LaTeX\ $\backslash${\tt{frac}}\{$\ldots$\}\{$\ldots$\}
and $\backslash${\tt{sqrt}}\{$\ldots$\} constructs~--
see for example (5.5)--(5.7).
\medskip

\noindent
A preliminary
version (rpb097tr) appeared as
Report CMA-R32-85, Centre for Mathematical Analysis,
Australian National University,
September 1985.
It is more detailed but does not include the
section on ``further refinements'' (\S9 above).
\medskip

\noindent
For developments up to mid-1997, see:

\begin{itemize}
\item[30.] R.~P.~Brent,
Factorization of the tenth Fermat number,
{\em Mathematics of Computation} 68~(January 1999), to appear (rpb161).
A preliminary version
({\em Factorization of the tenth and eleventh Fermat numbers},
Technical Report
TR-CS-96-02, Computer Sciences Laboratory, ANU, February 1996)
is also available in electronic form (rpb161tr).
\end{itemize}

\pagebreak[3]
\section*{Further remarks (added 3 December 1998)}

In the estimate (1.1), $T_1(p)$ is the arithmetic complexity.
The bit complexity is a factor $M(N)$ larger, where $M(N)$ is the
number of bit operations required to multiply integers mod~$N$.
As explained in~\S2,
we take one multiplication mod~$N$ as the basic unit of work.
In applications such as the factorization of large Fermat numbers,
the factor $M(N)$ is significant.

In \S4 the group operation on the elliptic curve
is written as multiplication, because of
the analogy with the Pollard ``$p-1$'' method.
Nowadays the group operation
is nearly always written as addition, see for example~[30].

At the end of \S8, we say that ``the product $mT$ is relatively insensitive
to the choice of~$m$''.  See [30, Table~3] for an indication of how
the choice of non-optimal $m$ changes the efficiency of the method.

\subsection*{Acknowledgement}

I am grateful to Paul Zimmermann for his comments which prompted
these ``further remarks''.

\end{document}